\documentstyle[11pt]{article}   
\newtheorem{thm}{Theorem}[section]

\newtheorem{corollary}[thm]{Corollary}
\newtheorem{proposition}[thm]{Proposition}
\newtheorem{definition}[thm]{Definition}
\newtheorem{example}[thm]{Example}
\newtheorem{remark}{Remark}

\newcommand{\qed}{\hfill\rule{2mm}{2mm}}  
\newenvironment{proof}{
\begin{trivlist}
\item[\hspace{\labelsep}{\bf\noindent Proof. }]
}{\qed\end{trivlist}}
\title{\Large \bf An umbral setting for cumulants and factorial moments}
\author{
 {\sc E. Di Nardo, D. Senato \footnote {Dipartimento di Matematica, Universit\`a degli
 Studi della Basilicata, Potenza, Italy. \{dinardo,senato\}@unibas.it}}}
 \date{\empty}

\begin{document}
\setlength{\baselineskip}{14pt}
\maketitle
\begin{abstract}
We provide an algebraic setting for cumulants and factorial moments
through the classical umbral calculus. Main tools are the compositional
inverse of the unity umbra, connected with the logarithmic power series,
and a new umbra here introduced, the singleton umbra. Various formulae
are given expressing cumulants, factorial moments and central moments by umbral
functions.
\\ \medskip \\
{\bf Keywords:\/} Umbral calculus, generating function, cumulant, factorial moment
\\
{\bf MSC-class:} 05A40, 60C05 (Primary)
62A01(Secondary)
\end{abstract}
\section{Introduction}
%
The purpose of this paper is mostly to show how the classical umbral
calculus gives a lithe algebraic setting in handling cumulants and
factorial moments. The classical umbral calculus consists
of a symbolic technique dealing with sequences of numbers $a_n$
indexed by nonnegative integers $n=0,1,2,3,\ldots,$ where
the subscripts are treated as if they were powers. This kind of
device was extensively used since the nineteenth century although
the mathematical community was sceptic of it, owing to its
lack of foundation. To the best of our knowledge, the method
was first proposed by Rev. John Blissard
in a series of papers as from $1861$ (cf. \cite{Dinardo}
for the full list of papers), nevertheless it is
impossible to put the credit of the original idea
down to him since the Blissard's calculus has its mathematical
source in symbolic differentiation. In the thirties, Bell \cite{Bell4} reviewed
the whole subject in several papers, restoring the purport of the Blissard's idea
and in \cite{Bell2} he tried to give a rigorous foundation of
the mystery at the ground of the umbral calculus but his attempt
did not have a hold. Indeed, in the first modern textbook of
combinatorics \cite{Riordan} Riordan often employed this
symbolic method without giving any formal justification.
It was first Gian-Carlo Rota to disclose the
\lq \lq umbral magic art\rq\rq$\,$ of shifting from $a^n$ to $a_n$
bringing to the light the underlying linear functional (cf. \cite{Rota}).
This idea led Rota and his collaborators to conceive
a beautiful theory (cf. \cite{Mullin} and \cite{VIII}) which
has originated a large variety of applications (see \cite{dibucchianico} for
a list of papers updated to $2000$). Some years later, Roman and Rota
\cite{Roman} gave rigorous form to the umbral tricks
in the setting of Hopf algebra (see also \cite{Joni}). But in $1994$ Rota himself wrote
(cf. \cite{SIAM}): {\small \lq \lq ... Although
the notation of Hopf algebra satisfied the most ardent advocate of
spic-and-span rigor, the translation of \lq \lq classical\rq\rq$\,$
umbral calculus into the newly found rigorous language made the
method altogether unwieldy and unmanageable. Not only was the eerie
feeling of witchcraft lost in the translation, but, after such a
translation, the use of calculus to simplify computation and sharpen our intuition
was lost by the wayside...\rq\rq$\,$}  Then, in the paper  \cite{SIAM} {\it The Classical
Umbral Calculus} (1994) Rota, together with Taylor, tries to restore the
feeling meant by the founders of the umbral calculus keeping new
notation both minimal and indispensable to avoid the misunderstanding of
the past. In this new setting, the basic device is to represent an unital sequence
of numbers by a symbol $\alpha,$ named umbra, i.e. to associate the sequence $1,a_1,
a_2, \ldots$ to the sequence $1,\alpha,\alpha^2, \ldots$ of
powers of $\alpha$ through an operator $E$ that looks like the
expectation of random variables (r.v.'s). This new way of dealing with
sequences of numbers has been applied to combinatorial and algebraic
subjects (cf. \cite{RotaTaylor}, \cite{Taylor} and \cite{Gessel}),
wavelet theory (cf. \cite{Shen1}) and difference equations
(cf. \cite{Taylor1}). Besides it has led to a nimble language
for r.v.'s theory, as showed in \cite{Shen} and \cite{Dinardo}.
\par
The present work is inspired by this last point of view. As a matter of fact,
an umbra looks as the framework of a random variable (r.v.)
with no reference to any probability space, someway getting
closer to statistical methods. However, the use of symbolic methods
in statistics is not a novelty: for instance Stuart and
Ord \cite{Stuart-Ord} resort to such a technique in handling moments about a point.
In addition in the umbral calculus, questions as convergence of series
are no matter, as showed hereafter dealing with cumulants.
\par
Among the sequences of numbers related to r.v.'s, cumulants play a central
role characterizing all r.v.'s occurring in the classical stochastic
processes. For instance, a r.v. having Poisson distribution of
parameter $x$ is the unique probability distribution
for which all its cumulants are equal to $x$. It seems
therefore that a r.v. is better described by its cumulants
than by its moments. Moreover, due to their properties of additivity and invariance under
translation, the cumulants are not necessarily connected with the moments of
any probability distribution. We can define cumulants $\kappa_j$
of any sequence $a_n, \, n = 1, 2, 3, ... $ by
$$\sum_{n=0}^{\infty} \frac{a_n t^n}{n!} = \exp \left\{
\sum_{j=1}^{\infty} \frac{\kappa_j t^j}{j!} \right\}$$
in disregard of questions of whether any series converges.
By this approach, many difficulties connected to the \lq\lq problem
of cumulants\rq\rq \, smooth out, where with \lq\lq problem of cumulants\rq\rq
\, we refer to characterizations of sequences that are cumulants of
some probability distributions. The simplest example is that the second
cumulant of a probability distribution must always be nonnegative,
and is zero only if all of the higher cumulants are zero.
Cumulants are subject to no such constraints when they are
analyzed by an algebraic point of view. What is more, in statistics
they do not play any dual role compared to factorial moments. Whereas the
algebraic setting here proposed comes to the light their close
relationship through an umbral analogy with the well known complementary
notions of compound and randomized Poisson r.v.'s (cf. \cite{FellerII})
\par
Umbral notations are introduced in Section 2 by means of
r.v.'s semantics. Our intention by this way is to make the reader
comfortable with the umbral system of calculation without require
any prior knowledge. We only skip some technical proofs
of formal matters on which the reader is referred to citations.
We also resume the theory of Bell umbrae, completely developed in \cite{Dinardo},
that not only gives the umbral counterpart of the family of Poisson r.v.'s
but it allows an umbral expression of the functional
composition of exponential power series. Section 3 is devoted to
a new umbra, named singleton umbra, playing a dual role compared to the Bell
umbra. Their relationship is encoded by the compositional inverse of the unity umbra.
The singleton umbra is the keystone of the umbral
presentation of cumulants and factorial moments.
\par
In the last two sections we give various umbral formulae for cumulants and factorial moments that
parallel those known in statistics but simplify the proofs as well as the forms.
This happens for instance for the equations expressing cumulants in terms of moments
(and vice-versa) and also for their recursive formulas. Inversion theorems
allowing to obtain an umbra from its cumulants or factorial moments are also stated.
\par
In $1929,$ Fisher \cite{Fisher} introduced the $k-$statistics as new symmetric functions
of the random sample. The aim of Fisher was to estimate the cumulants
without using the moment estimators. He used only combinatorial methods.
The $k-$statistics are related to the power sum symmetric functions whose
variables are the r.v.'s of the sample, but these expressions
are very unhandy. We believe that the umbral calculus may seek to simplify
the expression of the $k-$statistics (as well as the $h-$statistics for the central
moments) taking into account its combinatorial nature.
\section{Umbrae and random variables}
In the following, we resume terminology, notations and some basic
definitions of the classical umbral calculus, as it has been introduced by
Rota and Taylor in \cite{SIAM} and further developed in \cite{Dinardo}.
Fundamental is the idea of associating a sequence of numbers $1, a_2, a_3, \ldots$
to an indeterminate $\alpha$ which is said to represent the sequence.
This device is familiar in probability when $a_i$ represents
the $i-$th moment of a r.v. $X.$ In this case, the sequence $1, a_1, a_2, \ldots$
results from applying the expectation operator $E$  to the
sequence $1, X, X^2, \ldots$ consisting of powers of the r.v. $X.$
\par
More formal, an umbral calculus consists of the following data:
\begin{description}
\item[{\it a)}] a set $A=\{\alpha,\beta, \ldots \},$ called
the {\it alphabet}, whose elements are named {\it umbrae};
\item[{\it b)}] a commutative integral domain $R$ whose quotient field
is of characteristic zero;
\item[{\it c)}] a linear functional $E,$ called {\it evaluation},
defined on the polynomial ring $R[A]$ and taking values in $R$
such that
\begin{description}
\item[{\it i)}] $E[1]=1;$
\item[{\it ii)}] $E[\alpha^i \beta^j \cdots \gamma^k] =
E[\alpha^i]E[\beta^j] \cdots E[\gamma^k]$
for any set of distinct umbrae in $A$ and for $i,j,\ldots,k$
nonnegative integers ({\it uncorrelation property});
\end{description}
\item[{\it d)}] an element $\epsilon \in A,$ called {\it augmentation}
\cite{Roman}, such that $E[\epsilon^n] = \delta_{0,n},$
for any nonnegative integer $n,$ where
$$\delta_{i,j} = \left\{ \begin{array}{cc}
1 & \hbox{if $i=j$} \\
0 & \hbox{if $i \ne j$}
\end{array} \right. \,\, i,j \in N;$$
\item[{\it e)}] an element $u \in A,$ called {\it unity}
umbra \cite{Dinardo}, such that $E[u^n]=1,$ for any
nonnegative integer $n.$
\end{description}
A sequence $a_0=1,a_1,a_2, \ldots$ in $R$ is umbrally
represented by an umbra $\alpha$ when
$$E[\alpha^i]=a_i, \quad \hbox{for} \,\, i=0,1,2,\ldots.$$
The elements $a_i$ are called {\it moments} of the umbra
$\alpha$ on the analogy of r.v.'s theory. The umbra $\epsilon$
can be view as the r.v. which takes the value $0$ with
probability 1 and the umbra $u$ as the r.v. which takes
the value $1$ with probability 1. Note that
the uncorrelation property among umbrae parallels the analogue one
for r.v.'s as well as it is $E[\alpha^{n+k}] \ne E[\alpha^n] E[\alpha^k].$
Remark as this setting gets out of the well-known \lq\lq moment problem\rq\rq
\, for r.v.'s.   \par
\begin{example} {\it Bell umbra.} \label{exbell} \\
{\rm The {\it Bell umbra} $\beta$ is the umbra such that
$$E[(\beta)_n] = 1 \quad n=0,1,2,\ldots$$
where $(\beta)_0 =1$ and $(\beta)_n=\beta(\beta-1)\cdots(\beta-n+1)$
is the lower factorial. It results $E[\beta^n]= B_n$ where $B_n$
is the $n-$th Bell number (cf. \cite{Dinardo}), i.e. the number of the
partitions of a finite nonempty set with $n$ elements or the $n-$th coefficient
in the Taylor series expansion of the function $\exp(e^t-1).$
So $\beta$ is the umbral counterpart of the Poisson r.v. with parameter
$1.$}
\end{example}
We call {\it factorial moments} of an umbra $\alpha$ the
elements
$$a_{(n)} = \left\{\begin{array}{ll}
1, & n=0 \\
E[(\alpha)_n], & n>0
\end{array}\right.$$
where $(\alpha)_n=\alpha(\alpha-1)\cdots(\alpha-n+1)$ is the lower
factorial. So the definition of $\beta$ in example \ref{exbell}
could be reformulated as follows: the Bell scalar umbra is the
umbra whose factorial moments are $b_{(n)}=1$ for any nonnegative integer $n.$
\subsection{Similar umbrae and dot-product}
The notion of similarity among umbrae comes in handy in order
to manipulate sequences such
\begin{equation}
\sum_{i=0}^n \left( \begin{array}{c} n \\ i \end{array} \right)
a_{i}a_{n-i}, \quad \hbox{$n \in N$}
\label{(eq:1)}
\end{equation}
as moments of umbrae. The sequence (\ref{(eq:1)}) cannot be represented
by using only the umbra $\alpha$ with moments $a_0=1,a_1,a_2, \ldots.$
Indeed, being $\alpha$ correlated to itself, the product $a_i a_{n-i}$ cannot be
written as $E[\alpha^i \alpha^{n-i}].$ So we need two distinct
umbrae having the same sequence of moments, as it happens
for similar r.v.'s. Therefore, if we choose an umbra $\alpha^{\prime}$
uncorrelated with $\alpha$ but with the same sequence of moments,
it is
\begin{equation}
\sum_{i=0}^n \left( \begin{array}{c} n \\ i \end{array} \right)
a_{i}a_{n-i} = E\left[\sum_{i=0}^n \left( \begin{array}{c} n \\ i \end{array}
\right)
\alpha^i (\alpha')^{n-i}\right]= E[(\alpha + \alpha^{\prime})^n].
\label{(forment)}
\end{equation}
Then the sequence (\ref{(eq:1)}) represents the moments of the
umbra $(\alpha+\alpha^{\prime}).$ A way to formalize this matter is
to define two equivalence relations among umbrae.
\par
Two umbrae $\alpha$ and $\gamma$ are {\it umbrally equivalent} when
$$E[\alpha]=E[\gamma],$$
in symbols $\alpha \simeq \beta.$
They are {\it similar} when
$$\alpha^n \simeq \gamma^n, \quad n=0,1,2,\ldots$$
in symbols $\alpha \equiv \gamma.$ We note that equality
implies similarity which implies umbral equivalence. The
converses are false. Then, we shall denote by the symbol $n.\alpha$ the
{\it dot-product} of $n$ and $\alpha,$ an auxiliary umbra (cf. \cite{SIAM})
similar to the sum $\alpha^{\prime}+\alpha^{\prime
\prime}+ \ldots + \alpha^{\prime\prime\prime}$ where $\alpha^{\prime},
\alpha^{\prime \prime},\ldots, \alpha^{\prime\prime \prime}$
are a set of $n$ distinct umbrae each similar to the
umbra $\alpha.$ So the sequence in (\ref{(forment)}) is umbrally represented by the
umbra $2.\alpha$. We assume that $0.\alpha$ is an umbra similar to the augmentation $\epsilon.$
\par
We shall hereafter consider the dot product of $n$ and $\alpha$ as an umbra
if we {\it saturate} the alphabet $A$ with sufficiently many umbrae similar to
any expression whatever. For a formal definition of a saturated umbral
calculus see \cite{SIAM}. It can be shown that saturated umbral calculi
exist and that every umbral calculus can be embedded in a saturated
umbral calculus.
\par
The following statements are easily to be proved:
\begin{proposition}
\label{(prop1)}
\begin{description}
\item[\it{(i)}] If $n.\alpha \equiv n.\beta$ for some integer $n \ne 0$
then $\alpha \equiv \beta;$
\item[\it{(ii)}] if $c \in R$ then $n.(c \alpha) \equiv c (n.\alpha)$
for any nonnegative integer $n;$
\item[\it{(iii)}] $n.(m.\alpha) \equiv (nm).\alpha \equiv m.(n.\alpha)$
for any two nonnegative integers $n,m;$
\item[\it{(iv)}] $(n+m).\alpha \equiv n.\alpha + m.\alpha^{\prime}$
for any two nonnegative integers $n,m$ and any two distinct umbrae $\alpha \equiv \alpha^{\prime};$
\item[\it{(v)}] $(n.\alpha + n.\beta)\equiv n.(\alpha+\beta)$
for any nonnegative integer $n$ and any two distinct umbrae $\alpha$ and $\beta.$
\end{description}
\end{proposition}
Two umbrae $\alpha$ and $\gamma$ are said to be {\it inverse} to each other when
$\alpha+\gamma \equiv \varepsilon.$ We denote the inverse of the umbra $\alpha$ by
$-1.\alpha^{\prime},$ with $\alpha \equiv \alpha^{\prime}.$ Recall that,
in dealing with a saturated umbral calculus, the inverse of an umbra is not
unique, but any two inverse umbrae of the same umbra are similar.
\begin{example} {\it Uniform umbra.} \\
{\rm The {\it Bernoulli umbra} (cf. \cite{SIAM}) represents the sequence of
Bernoulli numbers $B_n$ such that
$$\sum_{k \geq 0} \left( \begin{array}{c}
n \\ k \end{array} \right) B_k = B_n.$$
The inverse of the Bernoulli umbra is the umbral counterpart of
the uniform r.v. over the interval $[0,1]$ (cf. \cite{Taylor})}.
\end{example}
\subsection{Generating functions}
The formal power series in $R[A][[t]]$
\begin{equation}
u + \sum_{n \geq 1} \alpha^n \frac{t^n}{n!}
\label{(gf)}
\end{equation}
is the {\it generating function} (g.f.) of the umbra $\alpha,$
and it is denoted by $e^{\alpha t}.$
The notion of umbrally equivalence and similarity can be extended
coefficientwise to formal power series $R[A][[t]]$
(see \cite{Taylor1} for a formal construction).
So it results
$$\alpha \equiv \beta \Leftrightarrow e^{\alpha t} \simeq e^{\beta
t}.$$
Moreover, any exponential formal power series\footnote{Observe that
with this approach we disregard of questions of whether any series converges.
} in $R[[t]]$
$$f(t) = 1 + \sum_{n \geq 1} a_n \frac{t^n}{n!}$$
can be umbrally represented by a formal power series (\ref{(gf)}) in $R[A][[t]].$
In fact, if the sequence $1,a_1,a_2,\ldots$ is umbrally
represented by $\alpha$ then
$$f(t)=E[e^{\alpha t}] \quad \hbox{i.e.} \quad f(t) \simeq e^{\alpha t},$$
assuming that we naturally extend $E$ to be linear.
We will say that $f(t)$ is umbrally represented by $\alpha.$ Note that, from now on, when
there is no mistaking, we will just say that $f(t)$ is the g.f. of $\alpha.$
For example the g.f. of the augmentation umbra $\epsilon$ is $1$ as well
as the g.f. of the unity umbra $u$ is $e^x.$
\par
Getting back to a r.v. $X,$ recall that when $E[\exp(t X)]$ is a convergent
function $f(t),$ it admits an exponential expansion in terms
of the moments which are completely determined by the related distribution function
(and vice-versa). In this case the moment generating function (m.g.f.)
encodes all the information of $X$ and the notion of similarity
among r.v.'s corresponds to that of umbrae.
\par
The first advantage of the umbral notation introduced for g.f.'s
is the representation of operations among g.f.'s with
operations among umbrae. For example the multiplication
among exponential g.f.'s is umbrally represented by a
summation of the corresponding umbrae:
\begin{equation}
g(t)f(t) \simeq e^{(\alpha+\gamma)t} \quad \hbox{with} \quad f(t) \simeq e^{\alpha
t}, \, g(t) \simeq e^{\gamma t}.
\label{(summation)}
\end{equation}
Via (\ref{(summation)}), the g.f. of $n.\alpha$ is $f(t)^n.$ If $\alpha$ is an umbra with g.f.
$f(t),$ the inverse umbra $-1.\alpha^{\prime}$ has g.f. $[f(t)]^{-1}.$
The summation among exponential g.f.'s is umbrally represented by a disjoint sum
of umbrae. The {\it disjoint sum} (respectively {\it disjoint difference}) of $\alpha$ and
$\gamma$ is the umbra $\eta$ (respectively $\iota$) with moments $$\eta^n
\simeq \left\{\begin{array}{ll}
u, & n=0 \\
\alpha^n + \gamma^n, & n>0
\end{array}\right. \quad \left( \hbox{respectively} \quad \iota^n \simeq \left\{\begin{array}{ll}
u, & n=0 \\
\alpha^n - \gamma^n, & n>0
\end{array}\right. \right), $$
in symbols $\eta \equiv \alpha \dot{+} \gamma$ (respectively $\iota \equiv \alpha \dot{-}
\gamma$). By the definition, it follows
$$f(t) \pm [g(t)-1] \simeq e^{(\alpha \dot{\pm} \gamma)t}.$$
\begin{example} {\it Unbiased estimators.} \label{exun} \\
{\rm Suppose to make the disjoint sum of $n$ times the umbra $\alpha.$
We will denote this umbra by $\dot{+}_{n}\alpha.$ Its
g.f. is $1+n[f(t)-1].$ The umbra $\dot{+}_{n}\alpha$
has the following probabilistic counterpart. Let $\{X_i\}_{i=1}^n$
be a random sample of independent and identically distributed (i.i.d.)
r.v's. As it is well-known  the power sum symmetric functions
$$S_r=\sum_{i=1}^n X_i^r$$
gives the unbiased estimators $S_r/n$ of the moments of $X_i.$
But $E[(\dot{+}_{n}\alpha)^r]=n \, a_r,$ hence the umbral
corresponding of the power sum symmetric functions sequence $S_r$ is
the umbra $\dot{+}_{n}\alpha.$}
\end{example}
\subsection{Auxiliary umbrae}
In the following, suppose $\alpha$ an umbra with g.f. $f(t)$
and $\gamma$ an umbra with g.f. $g(t).$ The introduction of
the g.f. device leads to the definition of new auxiliary umbrae
useful for the development of the system of calculation.
For this purpose, we should replace $R$ with whatever polynomial
ring having coefficients in $R$ and a number of indeterminates
according to necessity. In this paper, we deal with $R[x,y].$ This allows
to define the dot-product of $x$ and $\alpha$ via g.f., i.e. $x.\alpha$
is the auxiliary umbra having generating function
$$e^{(x.\alpha)} \simeq f(t)^x.$$
The Proposition \ref{(prop1)} still holds replacing $n$ with $x$ and $m$ with $y.$
Then, an umbra is said to be {\it scalar} if the moments are elements of $R$ while
it is said to be {\it polynomial} if the moments are polynomials.
\begin{example} {\it Bell polynomial umbra.} \\
{\rm The {\it Bell polynomial umbra} $\phi$ is the umbra having
factorial moments equal to $x^n$ (cf. \cite{Dinardo}). This umbra has
g.f. $\exp[x(e^t-1)]$ so that $\phi \equiv x.\beta,$
where $\beta$ is the Bell umbra. It turns out that the Bell polynomial umbra  $x.\beta$
is the umbral counterpart of the Poisson r.v. with parameter $x.$}
\end{example}
\begin{example}{\it Moments about a point.} \\
{\rm The moments $E[(X-a)^n]$ about a point $a \in {\bf R}$ of a r.v. $X$
are easily represented by umbrae through the following definition: the
umbra $\alpha^a$ having moments about a point $a \in R$
is defined as
\begin{equation}
\alpha^a \equiv \alpha - a.u.
\label{(cenmom)}
\end{equation}
If $a, b \in R$ and $b-a=c,$ then
$$\alpha^a \equiv \alpha-(b+c).u \equiv \alpha^b+c.u,$$
is the umbral version of the equations giving the moments
about $a$ in terms of the moments about $b$ (cf.
\cite{Stuart-Ord} for another symbolic expression).}
\end{example}
The dot-product $\gamma.\alpha$ of two umbrae is the auxiliary umbra having g.f.
$$e^{(\gamma.\alpha)t} \simeq  [f(t)]^{\gamma}
\simeq e^{\gamma \log f(t)} \simeq g\left[ \log f(t) \right].$$
The moments of the dot-product $\gamma.\alpha$ are (cf. \cite{Dinardo})
\begin{equation}
E[(\gamma.\alpha)^n] = \sum_{i=0}^n g_{(i)} B_{n,i}(a_1,a_2,\ldots, a_{n-i+1})
\quad n=0,1,2,\ldots
\label{(gr:2bis)}
\end{equation}
where $g_{(i)}$ are the factorial moments of the umbra $\gamma,$
$B_{n,i}$ are the (partial) Bell exponential polynomials (cf. \cite{Riordan})
and $a_i$ are the moments of the umbra $\alpha.$ Observe that
$E[\gamma.\alpha]= g_1 \, a_1 = E[\gamma]\,E[\alpha.]$
The following properties hold (cf. \cite{Dinardo}):
\begin{proposition} \label{prop2}
\begin{description}
\item[{\it a)}] if $\eta.\alpha \equiv \eta.\gamma$ then
$\alpha \equiv \gamma;$
\item[{\it b)}] if $c \in R$ then $\eta.(c \alpha) \equiv c (\eta.\alpha)$
for any two distinct umbrae $\alpha$ and $\eta;$
\item[{\it c)}] if $\gamma \equiv \gamma^{\prime}$
then $(\alpha+\eta).\gamma \equiv \alpha.\gamma + \eta.\gamma^{\prime};$
\item[{\it d)}] $\eta.(\gamma.\alpha) \equiv (\eta.\gamma).\alpha.$
\end{description}
\end{proposition}
Observe that from property {\it b)} it follows
\begin{equation}
\alpha.x \equiv \alpha.(x u) \equiv x(\alpha.u) \equiv x\alpha.
\label{(6bis)}
\end{equation}
\begin{remark}
{\rm The auxiliary umbra $\gamma.\alpha$ is the umbral version of a random sum. Indeed the m.g.f.
$g[\log f(t)]$ corresponds to the r.v.
$S_N = X_1+X_2+\cdots+X_N$ where $N$ is a discrete r.v.
having m.g.f. $g(t)$ and $X_i$ are i.i.d. r.v.'s having
m.g.f. $f(t).$ The right-distributive property of the dot-product $\gamma.\alpha$ runs in parallel with
the probability theory because the random sum $S_{N+M}$ is
similar to  $S_{N} + S_{M},$ where $N$ and $M$ are independent discrete
r.v.'s. The left-distributive property of the dot-product $\gamma.\alpha$
does not hold as well as it happens in the r.v.'s theory. In fact, let $Z=X+Y$
be a r.v. with $X$ and $Y$  independent r.v.'s. As it is easy to verify, the random sum
$S_N = Z_1+Z_2+\cdots+Z_N,$ with $Z_i$ i.i.d. r.v.'s similar to $Z,$
is not similar to the r.v. $S^{X}_N+S^{Y}_N$ with
$S^{X}_N=X_1+X_2+\cdots+X_N,$ and $X_i$ i.i.d. r.v.'s similar to
$X$ and with $S^{Y}_N=Y_1+Y_2+\cdots+Y_N$ and $Y_i$ i.i.d.
r.v.'s similar to $Y.$}
\end{remark}
\begin{example} {\it Randomized Poisson r.v.} \label{exrp} \\
{\rm Let us consider the Bell polynomial umbra $x.\beta.$ If in the
place of $x$ we put a generic umbra $\alpha,$
we get the auxiliary umbra $\alpha.\beta$ whose factorial moments
are
$$(\alpha.\beta)_n \simeq \alpha^n \quad n=0,1,2,\ldots$$
and moments given by the exponential umbral polynomials (cf. \cite{Dinardo})
\begin{equation}
(\alpha.\beta)^n \simeq \Phi_n(\alpha) \simeq \sum_{i=0}^n S(n,i) \alpha^i \quad
n=0,1,2,\ldots.
\label{(expBell)}
\end{equation}
Its g.f. is $f[e^t-1].$ The umbra $\alpha.\beta$ represents a random sum
of independent Poisson r.v.'s with parameter $1$ indexed by an integer r.v. $Y,$
i.e. a randomized Poisson r.v. with parameter $Y.$}
\end{example}
As suggested in \cite{VIII}, there is a connection between
compound Poisson processes and polynomial sequence of binomial type, i.e.
sequence $\{p_n(x)\}$ of polynomials with degree $n$ satisfying the
identities
$$p_n(x+y) = \sum_{i=0}^n \left( \begin{array}{c}
n \\
i
\end{array}\right)p_i(x) p_{n-i}(y)$$
for any $n$ (cf. for instance \cite{Mullin}).
Two different approaches can be found in \cite{Cerasoli} and in
\cite{Stam1}. A natural device to make clear this connection
is the {\it $\alpha-$partition} umbra $\beta.\alpha,$ introduced in \cite{Dinardo}.
Its g.f. is $\exp[f(t)-1]$ and it suggests to interpret a partition umbra
as a compound Poisson r.v. with parameter $1.$ As well-known, a
compound Poisson r.v. with parameter $1$ is introduced as a random sum $S_N = X_1 + X_2 +
\cdots + X_N$ where $N$ has a Poisson distribution with parameter $1$. The umbra
$\beta.\alpha$ fits perfectly this probabilistic notion taking into consideration
that the Bell scalar umbra $\beta$ plays the role of a Poisson r.v. with parameter $1.$
What's more, since the Poisson r.v. with parameter $x$ is umbrally represented
by the Bell polynomial umbra $x.\beta,$ a compound Poisson r.v. with parameter
$x$ is represented by the {\it polynomial $\alpha-$partition} umbra
$x.\psi \equiv x.\beta.\alpha$ with g.f. $\exp[x(f(t)-1)].$ The name
\lq \lq partition umbra\rq\rq$\,$ has a probabilistic ground.
Indeed the parameter of a Poisson r.v. is usually denoted by $x=\lambda t,$
with $t$ representing a time interval, so that when this interval is
partitioned into non-overlapping ones, their contributions
are stochastic independent and add to $S_N.$
This last circumstance is umbrally expressed by the relation
\begin{equation}
(x+y).\beta.\alpha \equiv x.\beta.\alpha + y.\beta.\alpha
\label{(eq:somma)}
\end{equation}
giving the binomial property for the polynomial sequence
represented by $x.\beta.\alpha.$  In terms of g.f.'s, the
formula (\ref{(eq:somma)}) means that
\begin{equation}
h_{x+y}(t) = h_x(t)h_y(t)
\label{(eq:somma1)}
\end{equation}
where $h_x(t)$ is the g.f. of $x.\beta.\alpha.$
Viceversa every g.f. $h_x(t)$ satisfying the
equality (\ref{(eq:somma1)}) is the g.f. of
a polynomial $\alpha-$partition umbra. The $\alpha-$partition
umbra represents the sequence of {\it partition polynomials}
$Y_n = Y_n(a_1,a_2,\ldots,a_n)$ (or complete Bell exponential polynomials
\cite{Riordan}), i.e.
\begin{equation}
E[(\beta.\alpha)^n] = \sum_{i=0}^n B_{n,i}(a_1,a_2,\ldots,
a_{n-i+1})= Y_n(a_1,a_2,\ldots,a_n),
\label{(comp)}
\end{equation}
where $a_i$ are the moments of the umbra $\alpha.$ Moreover every
$\alpha-$partition umbra satisfies the relation
\begin{equation}
(\beta.\alpha)^{n} \simeq \alpha'
(\beta.\alpha+\alpha')^{n-1}
\quad \alpha \equiv \alpha', \, n=0,1,2,\ldots
\label{(prpart)}
\end{equation}
and conversely (see \cite{Dinardo} for the proof). The previous
property will allow an useful umbral characterization
of the cumulant umbra (see corollary \ref{corcum} in section 4.)
The umbra $\beta.\alpha$ plays a central role also in the umbral representation
of the composition of exponential g.f.'s. Indeed, the {\it composition umbra}
of $\alpha$ and $\gamma$ is the umbra $\tau \equiv \gamma.\beta.\alpha.$
The umbra $\tau$ has g.f. $g[f(t)-1]$ and moments
\begin{equation}
E[\tau^n]= \sum_{i=0}^n g_i B_{n,i}(a_1,a_2,\ldots,a_{n-i+1})
\label{(momcomp)}
\end{equation}
with $g_i$ and $a_i$ moments of the umbra $\gamma$
and $\alpha$ respectively. We denote by $\alpha^{<-1>}$
the compositional inverse of $\alpha,$ i.e.
the umbra having g.f. $f^{-1}(t)$ such that $f^{-1}[f(t)-1]=f[f^{-1}(t)-1]=1+t.$
For an intrinsic umbral expression of the compositional inverse
umbra see \cite{Dinardo}, where it is also stated
an umbral version of the Lagrange inversion formula.
\begin{example} {\it Randomized compound Poisson r.v.} \label{(ex4)} \\
{\rm As already underlined in example \ref{exrp}, the umbra $\gamma.\beta$
represents a randomized Poisson r.v. Hence it is natural to
look at the composition umbra as a compound randomized Poisson
r.v., i.e. a random sum indexed by a randomized Poisson r.v.
Moreover, being $(\gamma.\beta).\alpha
\equiv \gamma.(\beta.\alpha)$ (cf. statement {\it d)} of
Proposition \ref{prop2}), the previous relation allows to see this r.v. from another
side: the umbra $\gamma.(\beta.\alpha)$ generalizes the concept of a random
sum of i.i.d. compound Poisson r.v. with parameter $1$ indexed
by an integer r.v. $X,$ i.e. a randomized compound Poisson r.v.
with random parameter $X.$}
\end{example}
At the end, the symbol $\alpha^{.n}$ denotes
an auxiliary umbra similar to the product $\alpha^{\prime}
\alpha^{\prime \prime} \cdots \alpha^{\prime\prime\prime}$
where $\alpha^{\prime},\alpha^{\prime \prime},\ldots,\alpha^{\prime\prime
\prime}$ are a set of $n$ distinct umbrae each similar to the
umbra $\alpha.$ We assume that $\alpha^{.0}$ is an umbra similar
to the unity umbra $u.$ The moments of $\alpha^{.n}$ are:
\begin{equation}
E[(\alpha^{.n})^k] =E[(\alpha^{k})^{.n}] = a_k^n, \quad
k=0,1,2,\ldots,
\label{(eq:10)}
\end{equation}
i.e. the $n-$th power of the moments of the umbra $\alpha.$ Thanks to this notation
in \cite{Dinardo}, the umbral expression of the Bell exponential polynomials was
given as follows:
\begin{equation}
B_{n,i}(a_1,a_2,\ldots,a_{n-i+1}) \simeq \left( \begin{array}{c}
n \\
i
\end{array}\right) \alpha^{.i} (i.\overline{\alpha})^{n-i}
\label{(eq:32)}
\end{equation}
whenever $a_1 \ne 0$ and where $\overline{\alpha}$ is the umbra with moments
\begin{equation}
\displaystyle{E\left[\overline{\alpha}^{n} \right]=
\frac{a_{n+1}}{a_1 (n+1)}}, \quad n=1,2,\ldots.
\label{(momover)}
\end{equation}
\begin{example}{\it The central umbra.} \\
{\rm We call {\it central umbra} the umbra $\alpha^{a_1}$ having moments
about $a_1 = E[\alpha].$ From (\ref{(cenmom)}), the classical
relation between moments and central moments of a r.v. has
the following umbral expression:
$$(\alpha^{a_1})^n \simeq \sum_{k=0}^n \left( \begin{array}{c}
n \\
k
\end{array}\right) (-1)^{n-k} (\alpha^{\prime})^k \alpha^{.(n-k)}, \quad
\alpha \equiv \alpha^{\prime} \,\, n=1,2,\ldots$$
being $E[(a_1.u)^k]=a_1^k=E[\alpha^{.k}]$ from (\ref{(eq:10)}).}
\end{example}
%
\section{The singleton umbra}
The singleton umbra plays a dual role compared to the Bell
umbra, even if it has not a probabilistic counterpart.
Besides, the singleton umbra turns out to be an effective
symbolic tool in order to umbrally represent some
well-known r.v.'s as well as cumulants and factorial moments.
\begin{definition}[The singleton umbra] \label{def}
An umbra $\chi$ is said to be a {\it singleton} umbra
if
$$\chi^n \simeq \delta_{1,n} \quad n=1,2,\cdots.$$
\end{definition}
The g.f. of the singleton umbra $\chi$ is $1+t.$
\begin{example} {\it Gamma r.v.} \\
{\rm The m.g.f. of a Gamma r.v. with parameters $a$ and $c$ is
$$M(t) = \frac{1}{(1 - c t)^a}.$$
This g.f. is umbrally represented by the inverse of $-c(a.\chi)$ (see
{\it (ii)} of Proposition \ref{(prop1)} replacing $n$ by $a \in
R$).}
\end{example}
Table 1 lights up the duality between the singleton umbra $\chi$
and the Bell umbra $\beta$.
\par \smallskip
{\centering
\begin{tabular}{|c|c|}
\hline
Umbra & Generating function  \\ \hline
$\chi$ & $(1+t)=1+ \sum_{n=1}^{\infty} \left[
\sum_{k=1}^n s(n,k)\right]\frac{t^n}{n!}$ \\
$x.\chi$ &  $(1+t)^x = 1+ \sum_{n=1}^{\infty} \left[ \sum_{k=1}^n  s(n,k) x^k \right]
\frac{t^n}{n!}$ \\
$\alpha.\chi$ &  $f[\log(1+t)]$ \\
$\chi.\alpha$ & $1+\log[f(t)]$ \\
\hline
$\beta$ & $e^{e^t-1}=1+\sum_{n=1}^{\infty} \left[ \sum_{k=1}^n S(n,k)\right]
\frac{t^n}{n!}$ \\
$x.\beta$ & $e^{x(e^t-1)}=1+ \sum_{n=1}^{\infty} \left[ \sum_{k=1}^n S(n,k)
x^k \right] \frac{t^n}{n!}$ \\
$\alpha.\beta$ & $f(e^t-1)$ \\
$\beta.\alpha$ & $e^{f(t)-1}$ \\
\hline
\end{tabular}
\begin{center} Table 1. \end{center}}
\par
The connection between the singleton umbra $\chi$ and the Bell umbra $\beta$
is made clear in the following proposition.
\begin{proposition}
\label{(cumbell)}
Let  $\chi$ be the singleton umbra, $\beta$ the Bell umbra and $u^{<-1>}$
the compositional inverse of the unity umbra $u.$ It results
\begin{eqnarray}
& & \chi  \equiv u^{<-1>}.\beta \equiv \beta.{u^{<-1>}}, \label{(rel)} \\
& & \beta.\chi \equiv u \equiv \chi.\beta. \label{(rel1)}
\end{eqnarray}
\end{proposition}
\begin{proof} The g.f. of $u^{<-1>}.\beta.u$
is $1+t,$ being $u^{<-1>}$ and $u$ compositional inverses. So
equivalence (\ref{(rel)}) follows by property {\it a)} of
proposition \ref{prop2} being
$$u^{<-1>}.\beta.u \equiv \chi \equiv \chi.u$$
Equivalence (\ref{(rel1)}) follows via g.f.'s in Table 1.
\end{proof}
Distributive properties of the singleton umbra respect to the sum
and the disjoint sum of umbrae are given in the following.
\begin{proposition}
\label{(propcum)}
It results
\begin{eqnarray}
\chi.(\alpha+\gamma)  & \equiv &  \chi.\alpha \dot{+} \chi.\gamma \label{(21)}\\
(\alpha \dot{+} \gamma).\chi & \equiv & \alpha.\chi \dot{+} \gamma.\chi.
\label{(22)}
\end{eqnarray}
\end{proposition}
\begin{proof}
Let $f(t)$ be the g.f. of $\alpha$ and $g(t)$ the g.f. of $\gamma.$
Equivalence (\ref{(21)}) follows observing that the g.f. of
$\chi.(\alpha+\beta)$ is $1+\log[f(t)g(t)]
= 1 + \log[f(t)]+\log[g(t)],$ i.e. the g.f. of
$\chi.\alpha \dot{+} \chi.\beta.$ Equivalence (\ref{(22)}) follows observing that the g.f. of
$(\alpha \dot{+} \beta).\chi$ is $f[\log(1+t)]+g[\log(1+t)]-1,$
i.e. the g.f. of $\alpha.\chi \dot{+} \beta.\chi.$
\end{proof}
The notion of mixture of r.v.'s has an umbral counterpart
in the disjoint sum $\dot{+}.$ Indeed let $\{\alpha_i\}_{i=1}^n$ be
$n$ umbrae and $\left\{p_i \right\}_{i=1}^n \in R$ be $n$ weights such that
$$\sum_{i=1}^n p_i=1.$$
The mixture umbra $\gamma$ of $\{\alpha_i\}_{i=1}^n$ is the following
weighted disjoint sum of $\{\alpha_i\}_{i=1}^n$
\begin{equation}
\gamma \equiv \chi.p_1.\beta.\alpha_1 \dot{+} \chi.p_2.\beta.\alpha_2
\dot{+} \ldots \dot{+} \chi.p_n.\beta.\alpha_n
\label{(mixt)}
\end{equation}
where $\beta$ is the Bell umbra and $\chi$ is the singleton
umbra. From (\ref{(21)}) equivalence (\ref{(mixt)}) can be rewritten
as
$$\gamma \equiv \chi.(p_1.\beta.\alpha_1 + p_2.\beta.\alpha_2
+ \ldots + p_n.\beta.\alpha_n).
$$
Since the g.f. of $\sum_{i=1}^n p_i.\beta.\alpha_i$
is $\exp(\sum_{i=1}^n p_i[f_i(t)-1]),$ where $f_i(t)$ is
the g.f. of $\alpha_i,$ from Table 1 it follows that the g.f. of $\gamma$
is $ \sum_{i=1}^n p_i f_i(t).$
\par
\begin{example} {\it Bernoulli umbral r.v.} \label{(ex2)} \\
{\rm Let us consider the Bernoulli r.v. $X$ of parameter
$p.$ Its m.g.f. is $g(t)=q+p \, e^t$ with $q=1-p.$
The Bernoulli umbral r.v. is the mixture of the umbra}
$\varepsilon$ and the unity umbra $u:$
$$\xi \equiv \chi.q.\beta.\varepsilon \dot{+} \chi.p.\beta.u.$$
{\rm Recalling that $\chi.q.\beta.\varepsilon \equiv \varepsilon$ it
is}
$$\xi \equiv \chi.p.\beta.$$
{\rm Indeed it is}
$$E[e^{\xi \, t}] = 1+\log[e^{p(e^t-1)}] = q + p \, e^t.$$
\end{example}
\begin{example} {\it Binomial umbral r.v.} \label{(ex3)}\\
{\rm As it is well-known a binomial r.v. $Y$ with parameters $n
\in {\bf N}, \, p \in [0,1],$ is the sum of $n$ i.i.d. Bernoulli
r.v.'s having parameter $p.$ Then the binomial umbral r.v. is
$$n.\xi \equiv n.\chi.p.\beta.$$
The parallelism is evident if we recall that the m.g.f. of
the binomial r.v. $Y$ is $f(t)=(q+pe^t)^n.$}
\end{example}
\section{The cumulant umbra}
For a r.v. having moments $a_1,a_2,\ldots,a_n$ and cumulants
$\kappa_1,\kappa_2,\ldots,\kappa_n$ it is
\begin{equation}
a_n=\sum_{\pi} c_{\pi} \kappa_{\pi} \quad \hbox{and} \quad
\kappa_n=\sum_{\pi} d_{\pi} a_{\pi}
\label{(conn)}
\end{equation}
the sums here are taken over the partitions $\pi = [j_1^{m_1},j_2^{m_2},
\ldots,j_k^{m_k}]$ of the integer $n,$ and
\begin{eqnarray*}
c_{\pi} & = & \frac{n!}{(j_1!)^{m_1}(j_2!)^{m_2}\cdots(j_k!)^{m_k}}
\frac{1}{m_1!m_2!\cdots m_k!} \\
d_{\pi} & = & c_{\pi} (-1)^{\nu_{\pi}-1} (\nu_{\pi}-1)! \quad
\hbox{and} \quad \nu_{\pi} = m_1 + m_2 + \cdots + m_k \\
a_{\pi} & = & \prod_{j \in \pi} a_j \quad
\hbox{and} \quad \kappa_{\pi}  =  \prod_{j \in \pi} \kappa_j.
\end{eqnarray*}
In this section we show how the umbral calculus simplifies the
above expressions, as well as the recursive formulae which
give moments in terms of cumulants.
\par
Let $\alpha$ be an umbra with g.f. $f(t).$
\begin{definition} \label{defcum}
The cumulant of an umbra $\alpha$ is the umbra $\kappa_{\alpha}$
defined by
$$\kappa_{\alpha} \equiv \chi.\alpha$$
where $\chi$ is the singleton umbra.
\end{definition}
Definition \ref{defcum} gives the umbral version of the second
equality in (\ref{(conn)}). Moreover the first moment of the cumulant umbra
$\kappa_{\alpha}$ is $a_1,$ i.e. the first moment of the umbra $\alpha,$
being $E[\kappa_{\alpha}] = E[\chi] \, E[\alpha] = E[\alpha] = a_1.$
\begin{example} {\it Cumulant of the umbra $\varepsilon.$} \\
{\rm Since $\varepsilon \equiv \chi.\varepsilon,$ the umbra $\varepsilon$ is the
cumulant umbra of itself, i.e. $\kappa_{\varepsilon} \equiv
\varepsilon.$}
\end{example}
\begin{example} {\it Cumulant of the umbra $u.$} \\
{\rm Since $\chi \equiv \chi.u,$ the umbra $\chi$ is the cumulant umbra
of the umbra $u,$ i.e. $\kappa_u \equiv \chi.$}
\end{example}
\begin{example} {\it Cumulant of the Bell umbra.} \label{(ex1)}\\
{\rm Since $u \equiv \chi.\beta$ (see (\ref{(rel1)})), the umbra $u$ is the
cumulant umbra of the Bell umbra $\beta,$ i.e. $\kappa_{\beta} \equiv u.$
From example \ref{exbell}, the Poisson r.v. of parameter $1$ has cumulants equal to
$1.$}
\end{example}
\begin{proposition}
The cumulant umbra $\kappa_{\alpha}$ has g.f.
\begin{equation}
k(t) = 1 + \log[f(t)].
\label{(gfcum)}
\end{equation}
\end{proposition}
\begin{proof}
See Table 1.
\end{proof}
\begin{example}{\it Cumulant of the singleton umbra.} \label{exsi}\\
{\rm Since $1+\log(1+t)$ is the g.f. of the umbra $u^{<-1>},$
this umbra is the cumulant umbra of the umbra $\chi,$ i.e. $\kappa_{\chi}
\equiv u^{<-1>}.$}
\end{example}
\begin{example} {\it Cumulant of the Bernoulli umbral r.v.} \\
{\rm From example \ref{(ex2)}, the cumulant umbra of the Bernoulli umbral r.v.
is $\chi.(\chi.p.\beta) \equiv u^{<-1>}.p.\beta.$}
\end{example}
\begin{example} {\it Cumulant of the Binomial umbral r.v.} \\
{\rm From example \ref{(ex3)}, the cumulant umbra of the Binomial umbral r.v.
is $\chi.(n.\chi.p.\beta),$ i.e.
$$\chi.(\chi^{\prime}.p.\beta^{\prime}+\chi^{\prime \prime}.p.\beta^{\prime \prime}
+\cdots+\chi^{\prime \prime \prime}.p.\beta^{\prime \prime \prime}),$$
where $\chi^{\prime},\chi^{\prime \prime},\ldots,\chi^{\prime \prime \prime}$
are a set of $n$ distinct umbrae each similar to the singleton umbra $\chi$ as well as
$\beta^{\prime},\beta^{\prime \prime},\ldots,\beta^{\prime \prime \prime}$ are
a set of $n$ distinct umbrae each similar to the Bell umbra $\beta.$
From (\ref{(21)}) and recalling examples \ref{exun} and \ref{exsi},
it results
$$\chi.n.\chi.p.\beta \equiv \chi.\chi^{\prime}.p.\beta^{\prime} \dot{+}
\chi.\chi^{\prime \prime}.p.\beta^{\prime \prime} \dot{+} \cdots \dot{+}
\chi.\chi^{\prime \prime \prime}.p.\beta^{\prime \prime \prime} \equiv
\dot{+}_n u^{<-1>}.p.\beta.$$
This parallels the analogous result in probability theory.}
\end{example}
From (\ref{(gfcum)}), the moments of the cumulant umbra $\kappa_{\alpha}$
are
$$(k_a)_n=E[\kappa_{\alpha}^n]=\left[ \frac{d^n}{dt^n}\log\{f(t)\}\right]_{t=0}$$
that is equivalent to the definition of the $n-$th cumulant of
a r.v. $X$ having m.g.f. $f(t).$
\par
To state the explicit version of the second equality in (\ref{(conn)})
\begin{equation}
k_n=\sum_{i=1}^n (-1)^{i-1} (i-1)! B_{n,i}(a_1,a_2,\ldots,a_{n-i+1})
\label{(cummom)}
\end{equation}
giving cumulants in terms of moments, usually requires laborious
computations (cf. for example \cite{Niki}).
The umbral definition of cumulants allows
a simple proof of (\ref{(cummom)}). Indeed, being $\chi \equiv
u^{<-1>}.\beta,$ the cumulant umbra of $\alpha$ is
the umbral composition of $u^{<-1>}$ and $\alpha:$
$$\kappa_{\alpha} \equiv u^{<-1>}.\beta.\alpha$$
and then its moments are given by (\ref{(momcomp)}). Equality
(\ref{(cummom)}) follows recalling that the moments of
$u^{<-1>}$ are the coefficient of the exponential expansion
$$1+\log(1+t) = 1 + \sum_{i=1}^{\infty} (-1)^{i-1} (i-1)! \frac{t^i}{i!}.$$
Similarly, the three main algebraic properties of cumulants
can be easily recovered from next theorem.
\begin{thm}
It is
\begin{description}
\item[{\it a)}] (the additivity property)
\begin{equation}
\chi.(\alpha+\gamma)\equiv \chi.\alpha \dot{+} \chi.\gamma,
\label{(additivity)}
\end{equation}
i.e. the cumulant umbra of a sum of two umbrae is equal to the
disjoint sum of the two corresponding cumulant umbrae;
\item[{\it b)}] (the semi-invariance under traslation
property) for any $c \in R$
$$\chi.(\alpha+ c.u)\equiv \chi.\alpha \dot{+} \chi.c;$$
\item[{\it c)}] (the homogeneity property) for any $c \in R$
$$\chi.(c \, \alpha) \equiv c (\chi.\alpha).$$
\end{description}
\end{thm}
\begin{proof}
Property {\it a)} follows from (\ref{(21)}). Property {\it b)}
follows from (\ref{(additivity)}), setting
$\gamma \equiv c.u$ for any $c \in R.$ At the end, property {\it c)}
follows from {\it b)} of proposition \ref{prop2}.
\end{proof}
\begin{example} {\it Cumulant of the central umbra.} \\
{\rm The sequence of cumulants related to the central umbra $\alpha^{a_1}$
is the same of $\alpha$ excepting the first equal to $0.$
Indeed, by the additivity property of the cumulant umbra it is
$$\chi.(\alpha - a_1.u) \equiv \chi.\alpha \dot{-} \chi.a_1.$$
The results follows from (\ref{(6bis)}).}
\end{example}
The umbral version of the first equality in (\ref{(conn)}) is given in the following
theorem.
\begin{thm}[Inversion theorem]
Let  $\kappa_{\alpha}$ be the cumulant umbra of $\alpha,$ then
$$\alpha \equiv \beta.\kappa_{\alpha}$$
where $\beta$ is the Bell umbra.
\end{thm}
\begin{proof}
It is
$$\beta.\kappa_{\alpha} \equiv \beta.u^{<-1>}.\beta.\alpha \equiv \chi.\beta.\alpha
\equiv u.\alpha \equiv \alpha.$$
\end{proof}
The inversion theorem allows to calculate the moments of the umbra
$\alpha$ according to its cumulants. Recalling (\ref{(comp)})
it is
\begin{equation}
a_n=Y_n[(k_a)_1,(k_a)_2,\ldots,(k_a)_n]
\label{(momcum)}
\end{equation}
with $a_n$ the $n-$th moment of the umbra $\alpha$ and $(k_a)_n$
the $n-$th moment of the umbra $\kappa_{\alpha}.$ Equation
(\ref{(momcum)}) is the explicit version of the first equality
in (\ref{(conn)}).
\begin{remark}
{\rm The complete Bell polynomials in (\ref{(comp)}) are a
polynomial sequence of binomial type. Since from the
inversion theorem any umbra $\alpha$ could be seen
as the partition umbra of its cumulant $\kappa_{\alpha}$, it is
possible to prove a more general result: every polynomial
sequence of binomial type is completely determined by
its sequence of formal cumulants. Indeed, in \cite{Dinardo}
it is proved that any polynomial sequence of binomial
type represents the moments of a polynomial
umbra $x.\alpha$ and viceversa. So from the inversion
theorem any polynomial sequence of binomial type
represents the moments of a polynomial umbra
$x.\beta.\kappa_\alpha.$}
\end{remark}
The next corollary follows from (\ref{(prpart)}) and
from the inversion theorem.
\begin{corollary} \label{corcum}
If $\kappa_{\alpha}$ is the cumulant umbra of $\alpha,$ then
\begin{equation}
\alpha^n \simeq \kappa_{\alpha} (\kappa_{\alpha}+\alpha)^{n-1}
\label{(ric)}
\end{equation}
for any nonnegative integer $n.$
\end{corollary}
Equivalences (\ref{(ric)}) were assumed by Shen and Rota
in \cite{Shen} as definition of the cumulant umbra.
In terms of moments, equivalences (\ref{(ric)}) give
$$ a_n = \sum_{j=0}^{n-1} \left( \begin{array}{c}
n-1 \\
j \end{array} \right) a_j (k_a)_{n-j}$$
that is largely used in statistic framework \cite{Smith}.
\begin{example} {\it L\'evy process}. \\
{\rm Let $(X_t, t \geq 0)$ be a real-value L\'evy process, i.e. a
process starting from $0$ and with stationary and independent
increments. According to the L\'evy-Khintchine formula (cf.
\cite{FellerII}), if we assume that $X_t$ has a convergent
m.g.f. in some neighbourhood of $0,$ it is
\begin{equation}
E[e^{\theta X_t}] = e^{t k(\theta)}
\label{(levy)}
\end{equation}
where $k(\theta)$ is the cumulant g.f. of $X_1.$
The inversion theorem gives the umbral version of equation
(\ref{(levy)}):
$$t.\alpha \equiv t.\beta.\kappa_{\alpha}.$$}
\end{example}
\subsection{Cumulants of the Poisson r.v.'s}
From example \ref{(ex4)}, the umbra $\gamma.\beta.\alpha$
corresponds to a compound randomized Poisson r.v., i.e. a random sum $S_N=X_1+\cdots+X_N$
with $N$ a randomized Poisson r.v. of parameter the r.v. $Y.$
In particular $\alpha$ corresponds to $X$ and $\gamma$ corresponds to $Y.$
Since $\chi.(\gamma.\beta.\alpha)\equiv \kappa_{\gamma}.\beta.\alpha$
the cumulant umbra of the composition of $\alpha$ and $\gamma$
is the composition of $\alpha$ and $\kappa_{\gamma}.$
Then from (\ref{(momcomp)}), the cumulants of a compound randomized
Poisson r.v. are given by
\begin{equation}
\sum_{i=1}^n k_i B_{n,i}(a_1,a_2, \ldots,a_{n-i+1})
\label{(cumpoi)}
\end{equation}
where $a_i$ are the moments of the r.v. $X$ and $k_i$ are the
cumulants of the r.v. $Y.$ Now set $\gamma \equiv x.u$ in
$\gamma.\beta.\alpha.$ This means to consider a r.v. $Y$ such that
$P(Y=x)=1.$ Then, the random sum $S_N$ becomes a compound Poisson r.v. of
parameter $x$ corresponding to the polynomial $\alpha-$partition
umbra $x.\beta.\alpha,$ with $\alpha$ the umbral counterpart of $X.$
and cumulants
\begin{equation}
\sum_{i=1}^n k_i B_{n,i}(a_1,a_2, \ldots,a_{n-i+1}) \simeq x
a_n.
\label{(cumpoi1)}
\end{equation}
Indeed (\ref{(cumpoi1)}) follows from (\ref{(cumpoi)})
since the moments $k_i$ of $\chi.x$ are equal to $0,$ except the
first equal to $x.$ If $x=1$ the cumulant of $\alpha-$partition umbra is $\alpha$
so that the moments of $X$ are the cumulants of the corresponding compound
Poisson r.v.  Now, in $x.\beta.\alpha$ take
$\alpha \equiv u.$ From (\ref{(cumpoi1)}), the cumulants of the
Bell polynomial umbra $x.\beta$ are equals to $x$ as well as for the
Poisson r.v. of parameter $x.$
\par
At the end, in $\gamma.\beta.\alpha$ set $\alpha \equiv u.$ The cumulant
umbra of $\gamma.\beta$ is $\kappa_{\gamma}.\beta$
with $\kappa_{\gamma}$ the cumulant umbra of $\gamma.$ Its
probabilistic counterpart is a randomized Poisson r.v. of parameter the
r.v. $Y,$ corresponding to the umbra $\gamma.$ From (\ref{(expBell)})
the cumulants of a randomized Poisson r.v. of parameter the r.v. $Y$ are the
moments of $\kappa_{\gamma}.\beta,$ i.e.
$$\sum_{i=0}^n S(n,i) k_i$$
with $k_i$ the cumulants of the r.v. $Y.$
\section{The factorial umbra}
The factorial moments of a r.v. do not play a very prominent role in
statistics, but they provide very concise formulae for the
moments of some discrete distributions, like the binomial one.
\par
Let $\alpha$ be an umbra with g.f. $f(t).$
\begin{definition}
\label{(fact)}
An umbra $\varphi_{\alpha}$ is said to be an
$\alpha-$factorial umbra if
$$\varphi_{\alpha} \equiv \alpha.\chi$$
where $\chi$ is the singleton umbra.
\end{definition}
\begin{example} {\it $\varepsilon-$factorial umbra.} \\
{\rm Since $\varepsilon \equiv \varepsilon.\chi,$
the $\varepsilon-$factorial umbra is similar to the umbra $\varepsilon,$
i.e. $\varphi_{\varepsilon} \equiv \varepsilon.$}
\end{example}
\begin{example} {\it $u-$factorial umbra.} \\
{\rm Since $\chi \equiv u.\chi,$ the $u-$factorial umbra is similar to
the umbra $\chi,$ i.e. $\varphi_{u} \equiv \chi.$}
\end{example}
\begin{example} {\it $\beta-$factorial umbra.} \\
{\rm Since $u \equiv \beta.\chi$ from (\ref{(rel1)}), the
$\beta-$factorial umbra is similar to the unity umbra $u,$
i.e. $\varphi_{\beta} \equiv u.$
From example \ref{exbell}, the Poisson r.v. of parameter $1$
has factorial moments equal to $1.$}
\end{example}
\begin{example} {\it $\chi-$factorial umbra.} \\
{\rm From example \ref{exsi}, it is $\chi.\chi \equiv u^{<-1>}.$ The $\chi-$factorial umbra
turns out to be $u^{<-1>},$ i.e. $\varphi_{\chi} \equiv
u^{<-1>}.$}
\end{example}
\begin{proposition}
The $\alpha-$factorial umbra has g.f.
\begin{equation}
g(t)= f\left[ \log (1+t) \right].
\label{(gffact)}
\end{equation}
\end{proposition}
\begin{proof}
See Table 1.
\end{proof}
The $\alpha-$factorial umbra has moments equal to the factorial
moments of the umbra $\alpha,$ as the following proposition
shows.
\begin{proposition}
Let $\varphi_{\alpha}$ be an $\alpha-$factorial umbra. Then
$$\varphi_{\alpha}^n \simeq (\alpha)_n, \, n=0,1,2,\ldots.$$
\end{proposition}
\begin{proof}  By equation (\ref{(gr:2bis)}) and definition \ref{(fact)}  it is
\begin{equation}
E[(\varphi_{\alpha})^n] = E[(\alpha.\chi)^n] = \sum_{k=0}^n (a)_k B_{n,k}(\delta_
{1,1},\delta_{1,2},
\ldots, \delta_{1,n-k+1})
\label{(eq:3)}
\end{equation}
where $(a)_k$ are the factorial moments of the umbra $\alpha$
and $\delta_{1,i}$ are the moments of the umbra $\chi.$ By
(\ref{(eq:32)}) it results
$$B_{n,k}(\delta_{1,1},\delta_{1,2},\ldots, \delta_{1,n-k+1}) \simeq \left(
\begin{array}{c}
n \\
k
\end{array}\right) \chi^{.k} (k.\overline{\chi})^{n-k}.
$$
Since the umbra $\overline{\chi}$ has moments equal to $0$ and $\chi^{.k} \simeq 1,$
then
\begin{equation}
B_{n,k}(\delta_{1,1},\delta_{1,2},\ldots, \delta_{1,n-k+1})= \left\{ \begin
{array}{ll}
0, & \hbox{if $n > k$}\\
1, & \hbox{if $n = k.$}
\end{array}\right.
\label{(bnkchi)}
\end{equation}
Hence the equation (\ref{(eq:3)}) becomes $E[(\varphi_{\alpha})^n] = (a)_n.$
\end{proof}
\begin{example} {\it Factorial umbra of the central umbra.} \\
{\rm From property {\it c)} of Proposition \ref{prop2}, it is
$$\alpha^{a_1}.\chi \equiv (\alpha-a_1.u).\chi \equiv \alpha.\chi -
a_1.\chi' \equiv \varphi_{\alpha} - \varphi_{a_1.u},$$
with $\chi' \equiv \chi.$ Then the factorial
umbra of the central umbra $\alpha^{a_1}$ is the difference between the factorial
umbra of $\alpha$ and the factorial umbra of the umbra having
moments equal to $a_1.$ By (\ref{(gffact)}) its g.f. results
$f[\log(1+t)](1-t)^{a_1}.$}
\end{example}
\begin{example} {\it Factorial moments of the binomial r.v.} \\
{\rm Since the factorial moments characterize the binomial r.v., we
show how to evaluate them by umbral methods. As showed in example
\ref{(ex3)}, the umbral counterpart of the binomial r.v.
is $n.\chi.p.\beta.$ Due to (\ref{(rel1)}) and (\ref{(6bis)})
the corresponding factorial umbra is $n.(\chi.p.\beta).\chi \equiv n.\chi.p
\equiv p(n.\chi).$ Its g.f. is
$$g(t)=(1+t\,p)^n=\sum_{j=0}^n (n)_j p^j \frac{t^j}{j!}$$
and so the factorial moments are $(n)_j p^j.$ If $n=1,$ the factorial umbra is
$p \, \chi$ and from example \ref{(ex2)} the first factorial moment of the Bernoulli
r.v. is equal to $p$ while the others are equal to $0.$}
\end{example}
\begin{example} {\it Factorial umbra of the cumulant umbra.} \\
{\rm If $\kappa_{\alpha}$ is the cumulant umbra of $\alpha,$ then
$\kappa_{\alpha}.\chi$ is the factorial cumulant umbra of $\alpha,$  with
g.f. $1+\log[f(1+t)]$ by (\ref{(gffact)}).}
\end{example}
The following theorem allows to obtain the umbra $\alpha$ from
its factorial umbra $\varphi_{\alpha}.$
\begin{thm}[Inversion theorem]
Let  $\varphi_{\alpha}$ be the factorial umbra of $\alpha.$
It is
$$\alpha \equiv \varphi_{\alpha}.\beta$$
with $\beta$ the Bell umbra.
\end{thm}
\begin{proof} By the Proposition \ref{(cumbell)} it results
$$\varphi_{\alpha}.\chi \equiv \alpha.\chi.\beta \equiv \alpha.$$
\end{proof}
%
\subsection{Factorial moments of the Poisson r.v.'s}
%
Being $(\gamma.\beta.\alpha).\chi \equiv \gamma.\beta.(\varphi_{\alpha})$
the factorial umbra of the umbral composition $\gamma.\beta.\alpha$
is the umbral composition of $\gamma$ and the factorial umbra of $\alpha.$ From
(\ref{(momcomp)}) the compound randomized Poisson r.v. $S_N=X_1+\cdots+X_N$ with $N$
a Poisson r.v. with parameter the r.v. $Y$ has factorial moments
\begin{equation}
\sum_{k=1}^n g_k B_{n,k}[(\mu)_1,(\mu)_2, \ldots,(\mu)_{n-k+1}]
\label{(factpoi)}
\end{equation}
where $(\mu)_i$ are the factorial moments of the r.v. $X$ and $g_k$ are
the moments of the r.v. $Y$. Now setting $\gamma \equiv x.u$ in $\gamma.\beta.\alpha,$
we have $g_k=x^k.$ Then from (\ref{(factpoi)})
\begin{equation}
\sum_{k=1}^n x^k B_{n,k}[(\mu)_1,(\mu)_2, \ldots,(\mu)_{n-k+1}]
\label{(factpoi1)}
\end{equation}
are the factorial moments of a compound Poisson r.v. with parameter $x.$
Set $\alpha \equiv u$ in $x.\beta.\alpha.$ We have $(x.\beta.u).\chi \equiv x.\beta.\chi
\equiv x.u$ so that the factorial moments of $x.\beta.\alpha$
are equals to $x^n$ as well as for its probabilistic counterpart,
the Poisson r.v. with parameter $x.$
\par
At the end set $\alpha \equiv u$ in $\gamma.\beta.\alpha.$ We
have $(\gamma.\beta.u).\chi \equiv \gamma.\beta.\chi
\equiv \gamma$ so that the factorial moments of $\gamma.\beta$
are equals to the moments of $\gamma.$ Then a randomized Poisson
r.v. with parameter a r.v. $Y$ has factorial moments equal to
the moments of the r.v. $Y.$
\end{document}